\documentclass[a4paper,12pt]{article}
\usepackage{graphicx,amssymb}
\usepackage{amssymb,amsmath}
\usepackage{hyperref}
\usepackage{amsfonts}
\usepackage{amscd}
\usepackage{amsthm}
\usepackage{color}
\numberwithin{equation}{section}

\usepackage{amssymb,amsmath,amsthm,authblk}
\usepackage{enumitem}
\newtheorem{theorem}{Theorem}
\newtheorem{corollary}{Corollary}
\newtheorem{remark}{Remark}
\newtheorem{lemma}{Lemma}
\newtheorem{example}{Example}
\newtheorem{proposition}{Proposition}
\newtheorem{definition}{Definition}
\newtheorem{observation}{Observation}
\newtheorem{prob}{Problem} 

\begin{document}
\title{ \textbf{Spectrum of group vertex magic graphs}}
\baselineskip 16pt
\author{\large S. Balamoorthy $^\dagger$, N. Kamatchi$^\ddagger$ and S.V. Bharanedhar$^\dagger$\footnote{Corresponding Author. \\ E-Mail addresses: moorthybala545@gmail.com (S. Balamoorthy), \\ kamakrishna77@gmail.com (N. Kamatchi), bharanedhar3@gmail.com (S.V. Bharanedhar).} 
\\ {\small$^\dagger$Department of Mathematics, Central University of Tamil Nadu, Thiruvarur 610005, India}\\
{\small$^\ddagger$ Department of Mathematics, Kamaraj College of Engineering and Technology, Virudhunagar 625701, India}}
\date{\today}
\maketitle

\begin{abstract}
Let $G$ be a simple undirected graph and let $A$ be an additive Abelian group with identity 0. In this paper, we introduce the concept of group magic spectrum of a graph $G$ with respect to a given Abelian group $A$ and is defined as $spec(G,A)$:= $\{ \lambda: \lambda$ is a magic constant of some  $A$-vertex magic labeling $f\}$.  In their recent work, K. M. Sabeel et al. in Australas. J. Combin. 85(1) (2023), 49-60 proved a forbidden subgraph characterization for the group vertex magic graph. In this work, we present a new method which uses minimum number of vertices required for this graph. 
We obtain a necessary and sufficient condition for the spectrum of a graph $G$ to be a subgroup when $A=V_4$ or $\mathbb{Z}_p$, where $p$ is a prime number. Also we introduce the notion of reduced spectrum $redspec(G,A)$ and study the relation between $spec(G,A)$ and $redspec(G,A)$.

\medskip
\noindent {\bf  Keywords:} Group vertex magic; magic spectrum;  reduced graph; tensor product.\medskip 

\noindent {\bf  AMS Subject classification:} 05C25, 05C78, 05C76.
\end{abstract}

\section{Introduction}
Throughout this paper, we consider  finite, simple and connected graphs with vertex set $V(G)$ and edge set $E(G)$. For a vertex $v \in V(G)$, let $N_G(v)$ be the set of all vertices adjacent to $v$ in $G$ and $|N_G(v)| = deg_G(v)$. Let $V_4$ denote the famous Klein’s four group and it is well known that $V_4\cong\mathbb{Z}_2\times\mathbb{Z}_2$. 
Let $A$ be an additive Abelian group with identity 0. A mapping $ l : V(G) \to A \setminus \{0\}$ is said to be a $A$-vertex magic labeling of $G$ if there exists a $ \mu$ in $A$ such that  $w(v) = \sum_{u \in N_G(v)} l(u) = \mu $ for any vertex $v$ of $G$. If $G$ admits such a labeling, then it is called an $A$-vertex magic graph. If $G$ is $A$-vertex magic for any non-trivial Abelian group $A$, then $G$ is called a group vertex magic graph. Several works have been done on A-vertex magic graphs in \cite{Bala, Bala1, Bala2, Sabeel}.\\ 
Let $H$ be a graph with vertex set $\{v_1,v_2,\ldots,v_k\}$ and let $\mathcal{F}= \{G_1,G_2,\ldots,G_k \}$ be a family of graphs. The $H$-join operation of the graphs $G_1,G_2,\ldots,G_k$ is denoted as $G=H[G_1,G_2,\ldots,G_k]$, is obtained by replacing the vertex $v_i$ of $H$ by the graph $G_i$ for $1 \leq i \leq k$ and every vertex of $G_i$ is made adjacent with every vertex of $G_j$, whenever $v_i$ is adjacent to $v_j$ in $H$. The vertex set of $G$ is $V(G)= \displaystyle \bigcup_{i=1}^k V(G_i)$ and the edge set of $G$ is $E(G) =  \displaystyle \big(\bigcup_{i=1}^k E(G_i) \big) \cup \big(\bigcup_{v_iv_j \in E(H)} \{uv : u \in V(G_i), v \in V(G_j)\}\big)$ (see \cite{Vinoth} and \cite{Saravanan}). If $G'\cong G_i$, for $1\leq i\leq k$, then $H[G_1,G_2,\ldots, G_k]\cong H[G']$, is called the \emph{lexicographic product} of $H$ and $G'$.

It is always interesting to study graphs arising from algebraic structures. In this connection, one can see the recent survey articles \cite{Kumar} and \cite{Cameron}. In this direction, we introduce the concept of group magic spectrum of a graph $G$ with respect to a group $A$ as a set of group elements which comes as magic constant of a group magic labeling. So to explore this connection we introduce the following definition.  

\begin{definition}
If $G$ is an $A$-vertex magic graph for some Abelian group $A$, then the group magic spectrum with respect to the group $A$ is defined as 
$\{ \lambda: \lambda$ is a magic constant of some  $A$-vertex magic labeling f\} and is denoted by $spec(G,A).$
\end{definition}

That is, $spec(G,A)$ = $\{ \lambda: \lambda$ is a magic constant of some  $A$-vertex magic labeling $f$\}.  


From the above definition, it is clear that $spec(G,A)$ is a subset of $A$. We are interested in the following fundamental problem. 

\begin{prob}
For given a graph $G$ and an Abelian group $A$, when $spec(G,A)$ is a subgroup of $A$?
\end{prob}
    





\begin{observation}\cite{Bala1}\label{BP1115}
A graph $G$ is $\mathbb{Z}_2$ magic if and only if degree of every vertex in $G$ is of same parity.
\end{observation}
 
\begin{lemma}\cite{Bala2}\label{BP2011}
Let $A$ be an Abelian group with at least three elements. If $n \geq 2$ and $a \in A $, then there exists $a_1,a_2,\ldots,a_n$ in $A\setminus \{0\}$ such that $a = \sum_{i=1}^n a_i$. 
\end{lemma}

\begin{definition} The tensor product $G \otimes H$ of graphs $G$ and $H$ is a graph such that the vertex set of $G \otimes H$ is the product $V(G) \times V(H)$ and
vertices $(g,h)$ and $(g',h')$ are adjacent in $G \otimes H $ if and only if
$g$ is adjacent to $g'$ in $G$ and $h$ is adjacent to $h'$ in $H$.
\end{definition}

For basic graph-theoretic notion, we refer to Bondy and Murty \cite{bondy1976graph}. Let $R$ be a commutative ring with unity, we denote the multiplicative group of all units in $R$ by $U(R)$. For concepts in group theory, we refer to Herstein \cite{herstein}. 

\section{Main Results}
In this section, we present a few necessary conditions for $spec(G,A)$ to be a subgroup of a given group $A$.

 From the definition of tensor product one can observe the following. 

\begin{observation}
Let $G_1$  and $G_2$ be two simple graphs. The graphs $G_1$ and $G_2$ have odd degree vertices if and only if $G_1 \otimes G_2$ has odd degree vertices.  
\end{observation}
\begin{proof}
From the definition of tensor product, we have $deg_{G_1 \otimes G_2}(u,v) = deg_{G_1}(u)
\cdot deg_{G_2}(v)$, where $(u,v)\in V(G_1 \otimes G_2)$. Hence the result follows. \end{proof}

\begin{observation}\label{BP4009}
Let $G_i$ \, $(1 \leq i \leq k \, ; \, k \geq 2)$ be simple graphs and $G = \otimes_{i=1}^k G_i$. Then $G$ is Eulerian if and only if at least one of the ${G_i}'s$ is Eulerian.    
\end{observation}

\begin{lemma}
Let $A$ be an Abelian group, underlying a commutative ring $R$ and $a \in U(R)$. If $l$ is an $A$-vertex magic labeling of $G$ with magic constant $\mu$, then there exits an $A$-vertex magic labeling $l'$ of $G$ with magic constant $a\mu$.  
\end{lemma}
\begin{proof}
Assume that $l$ is an $A$-vertex magic labeling of $G$ with magic constant $\mu$. Define $l' : V(G) \to A\setminus \{0\}$ by $l'(v) = al(v)$. Clearly, $w(v)=a\mu$. Hence we get the required result.
\end{proof}

\begin{corollary}
Let $A$ be an Abelian group, underlying a commutative ring $R$. If $\mu \in spec(G,A)$, then $a\mu \in spec(G,A)$, where $a \in U(R)$.
\end{corollary}

\begin{theorem}\label{BP4008}
Let $G$ be an $A$-vertex magic graph. If $0\in spec(G,A)$, then $G$ has no pendant vertices.
\end{theorem}
\begin{proof}
Suppose $G$ has a
pendant vertex $v$ (say) and $N_G(v)= \{u\}$. Since $G$ is $A$-vertex magic, we have $w(v)=l(u) \neq 0$, where $l$ is any $A$-vertex magic labeling of $G$. Thus, by definition of $spec(G,A)$, $0$  is not an element of $spec(G,A)$.
\end{proof}

\begin{observation}\label{BP4013}
Let $G$ be an $A$-vertex magic graph. If $G$ has a pendant vertex, then $spec(G,A)$ is not a subgroup of $A$.     
\end{observation}

\begin{remark}
For no tree $T$, $spec(T,A)$ is a subgroup of $A$.
\end{remark}

\begin{theorem}\label{BP4010}
Let $G$ be a $\mathbb{Z}_2$-vertex magic graph. Then the $spec(G,\mathbb{Z}_2)$ is a subgroup of $\mathbb{Z}_2$ if and only if $G$ is Eulerian. In this case, $spec(G, \mathbb{Z}_2) = \{0\}.$
\end{theorem}
\begin{proof}
Assume that $spec(G,\mathbb{Z}_2)$ is a subgroup of $ \mathbb{Z}_2$. Since $G$ is a $\mathbb{Z}_2$-vertex magic graph, by Observation \ref{BP1115} and by group definition, we get the required result. Converse part is straightforward. 
\end{proof}


\begin{corollary}
 There is no graph $G,$ such that $spec(G,{\mathbb{Z}_2})= \mathbb{Z}_2$. 
\end{corollary}

\begin{theorem}
 Let $G_i$ \, $(1 \leq i \leq k \, ; \, k \geq 2)$ be simple graphs and let $G = \otimes_{i=1}^k G_i$. Then $spec(G,\mathbb{Z}_2)$ is subgroup of $\mathbb{Z}_2$ if and only if $spec(G_i,\mathbb{Z}_2)$  is subgroup of $\mathbb{Z}_2$, for some $i$.     
\end{theorem}
\begin{proof}
By Theorem \ref{BP4010} and Observation \ref{BP4009}, we get the required result.   
\end{proof}

\begin{theorem}\label{BP4015}
Let $G$ be a $V_4$-vertex magic graph. Then the $spec(G,V_4)$ is a subgroup of $V_4$ if and only if $0 \in spec(G,V_4)$. 
\end{theorem}
\begin{proof}
Clearly, it is enough to prove the sufficiency part. If $spec(G,V_4)=\{0\},$ then there is nothing to prove. Suppose $spec(G,V_4) \neq \{0\}$, then there exists a labeling $l$ with a non-zero magic constant. Assume that the magic constant of $l$ is $a$. 
Let $v \in V(G)$. Then $w(v) = r_1a+r_2b+r_3c,$ where $r_i \in \{0,1\},$ for all $i=1,2,3$. If $r_1a \neq 0,$ then $r_2b+r_3c =0$ and also $r_1a=a.$ If $r_1a =0,$ then $r_2b+r_3c=b + c=a$. We replace the labels $a,b,c$ in $l$ by $b,c,a$. Now, $w(v)= r_1b+r_2c+r_3a=b$. By a similar argument, we get $w(v)=c$. 
\end{proof}

\begin{corollary}
There is no graph $G,$ whose $spec(G,V_4)$ is a proper subgroup of $V_4$. 
\end{corollary}

The following theorem shows that $spec(G,A)$ is symmetric about the origin.
\begin{theorem}\label{BP4017}
If $a \in spec(G,A)$, then $-a \in spec(G,A)$ .
\end{theorem}
\begin{proof} 
Assume that $l$ is a $A$-vertex magic labeling of $G$ with magic constant $a$. Define $l' : V \to A \setminus \{0\}$ by $l'(v)= -l(v),$ for all $v \in V(G)$. Now, 
\begin{align*}
w(v) &= \sum_{u \in N_G(v)} l'(u)\\
&= -\sum_{u \in N_G(v)} l(u)\\
&= - a.
\end{align*}
Since $v$ is arbitrary, therefore, $w(v)= -a,$ for all $v \in V(G)$. 
\end{proof}
\begin{remark}
    Thus by the above theorem to check spectrum is a subgroup it is enough to check it is closed. 
\end{remark}

\begin{theorem}\label{BP4014}
Let $p$ be a prime number and $a \in \mathbb{Z}_p \setminus \{0\}$. If $a \in spec(G,\mathbb{Z}_p)$, then $\mathbb{Z}_p\setminus\{0\}\subset spec(G,\mathbb{Z}_p)$.
\end{theorem}
\begin{proof}
Let $a$ be a non-zero element such that $a \in spec(G, \mathbb{Z}_p)$. Then there exists a labeling $l$ such that the label of every vertex is a multiple of $a$ and $w(v)=a$, for all $v \in V(G)$. Let $v \in V(G)$. Since $a$ generates $\mathbb{Z}_p$, we have $$w(v) = r_1a+r_22a+\cdots+r_iia+\cdots+r_{p-1}(p-1)a, \, 0 \leq r_j \leq deg_G(v),$$ for all $j \in \{1,2,\ldots,p-1\}$.\\ Hence,
$w(v) = 
(r_1+2r_2+\cdots+ir_i+\cdots+(p-1)r_{p-1})a=a$.\\
Therefore, $(r_1+2r_2+\cdots+ir_i+\cdots+(p-1)r_{p-1}) \equiv 1({\rm mod} \, p)$. \\Hence $(r_1+2r_2+\cdots+ir_i+\cdots+(p-1)r_{p-1})ka =ka$, where $0 < k \leq p-1$.
Now, define $l': V \to \mathbb{Z}_p \setminus \{0\}$ by $l'(v_i)=kl(v_i)({\rm mod}\, p)$. Hence we have $ka \in spec(G,\mathbb{Z}_p)$.
\end{proof}

\begin{corollary}
Let $G$ be a $\mathbb{Z}_p$-vertex magic graph. Then the $spec(G,\mathbb{Z}_p)$ is a subgroup of $\mathbb{Z}_p$ if and only if $0 \in spec(G,\mathbb{Z}_p)$.
\end{corollary}

\noindent{From Theorem \ref{BP4015} and \ref{BP4014}, we have the following result}.

\begin{theorem}\label{BP4016}
Let $A = V_4$ or $\mathbb{Z}_p$ and $G$ be an $A$-vertex magic graph. The $spec(G,A)$ is a subgroup of $A$ if and only if $0 \in spec(G,A)$.
\end{theorem}

The following Proposition is an interesting result in group theory whose proof is trivial.

\begin{proposition}\label{BP4001}
Let $p_i$'s are distinct prime numbers  and $n=p_1^{\alpha_1}p_2^{\alpha_2}\ldots p_m^{\alpha_m}, \alpha_i > 0$ for all $i=1,2,\ldots,m$. The group $A$ has an element a such that $o(a)|n$ if and only if $A$ has a subgroup isomorphic to $\mathbb{Z}_{p_i}$ for some $i=1,2,\ldots,m$.
\end{proposition}

\noindent{Since the existence of $A$-vertex magicness of the graphs $C_n$ and complete $k$-partite graphs are studied in \cite{Bala1}, it is interesting to study their spectrum.}


\begin{theorem}\label{BP4007}
$Spec(C_n,A)$ is a subgroup of $A$ if and only if one of the following holds.
\begin{enumerate}
\item $n \equiv 0 (\textrm{mod}\ 4)$.
\item $n \not\equiv 0 (\textrm{mod}\ 4)$ and  $A$ has a subgroup isomorphic to $\mathbb{Z}_2$.
\end{enumerate}  
\end{theorem}
\begin{proof}
Let $V(C_n)= \{v_1,v_2,\ldots,v_n\}$. Since $C_n$ is $A$-vertex magic, a simple computation shows that,
\begin{equation}\label{BP4equ1}
l(v_i)=l(v_{(i+4)(\textrm{mod}\ n)}).
\end{equation} 
Assume that $spec(C_n,A)$ is a subgroup of $A$. Suppose $n \equiv 0 (\textrm{mod}\ 4)$, then there is nothing to prove. Suppose $n \not\equiv 0 (\textrm{mod}\ 4)$, then we have the following two cases.
\item[\bf Case 1.] $n$ is odd. \\ In this case by Equation \ref{BP4equ1}, we have $l(v_i)=l(v_j)$ for all $i,j \in \{1,2,\ldots,n\}$. Assume that $l(v_1)=a$, then $w(v_i)=2a$ for all $i$. Therefore $spec(C_n,A)= \{2a : a\in A\setminus \{0\}\}$. Since $spec(C_n,A)$ is subgroup of $A$, then $0 \in spec(C_n,A)$, which implies $A$ has an order $2$ element. 
\item[\bf Case 2.] $n \not\equiv 0(\textrm{mod}\ 4)$ and $n$ is even. In this case by Equation \ref{BP4equ1}, we have $l(v_1)=l(v_3)=l(v_5)=\cdots=l(v_{(n-1)})$ and $l(v_2)=l(v_4)=l(v_6)=\cdots=l(v_n)$. Assume that $l(v_1)=a$ and $l(v_2)=b$. Then \[ w(v_i)= \begin{cases} \text{$2a$} &\quad\text{if $i$ is even}\\ \text{$2b$} &\quad\text{if $i$ is odd.}\\ \end{cases} \] Therefore, the $spec(C_n,A)= \{2a : a \in A\setminus \{0\}\}$. Since $spec(C_n,A)$ is a subgroup of $A$, then $A$ has an order $2$ element.\\
Conversely, first let us assume that $n \not\equiv 0 (\textrm{mod}\ 4)$ and $A$ has a subgroup isomorphic to $\mathbb{Z}_2$. Since $A$ has an element of order $2$, which implies $ 0 \in spec(C_n,A)$. Let $ a,b \in spec(C_n,A)$, where $a = 2a_1, b = 2b_1$ and $a_1,b_1 \in A\setminus \{0\}$. Then $a-b = 2(a_1-b_1) \in spec(C_n,A)$. Next, assume that $n=4m$ 
and  $l(v_1)=a,l(v_2)=b,l(v_3)=c$ and $l(v_4)=d$. Since $w(v_2)=w(v_3)$, which implies $a+c =b+d$. Then the $spec(C_n,A) = \{a+b : a,b \in A\setminus \{0\} \}$. Let $c_1,c_2 \in spec(C_n,A)$, where $c_1 = a_1+b_1,c_2 = a_2+b_2$. Now, $c_1-c_2 = a_1-a_2 + b_1-b_2 \in spec(C_n,A)$. In this case $spec(C_n, \mathbb{Z}_2) = \{0\}$, otherwise $spec(C_n,A) = A$.
\end{proof}

\begin{proposition}\label{BP4002}
A complete $k$-partite graph $G$ is $A$-vertex magic if and only if sum of all labels of vertices in each partite set are equal.
\end{proposition}
\begin{proof}
Let $V_1,V_2,\ldots,V_k$ be a partition of $V(G)$ with $|V_j|=n_j$, for $j = 1,2,\ldots k$. Let $v_i^j \in V_j$. Since $G$ is $A$-vertex magic $w(v_i^1)=w(v_i^j)$, which implies $\sum_{i=1}^{n_j} l(v_i^j) = \sum_{i=1}^{n_1} l(v_i^1)$ for all $j$. Therefore, sum of all labels of vertices in each partite are equal. The proof of converse part is trivial.
\end{proof}

\begin{theorem}\label{BP4005}
Let $G$ be a complete $k$-partite graph with at least one partite size  equal to $1$, where $|A|>2$ and $k>2$. Then spec$(G,A)$ is subgroup of $A$ if and only if $A$ has a subgroup isomorphic to $\mathbb{Z}_p$, where $p$ is a prime number and $p| k-1$. \end{theorem}
\begin{proof}
Let $V_1,V_2,\ldots,V_k$ be a partition of $V(G)$ with $|V_j|=n_j$, for $j = 1,2,\ldots k$. Let $v_i^j \in V_j$.  As $G$ has at least one partite size is 1 and by Proposition \ref{BP4002}, we get the sum of all labels of vertices in each partite are equal to $a$ for some $a \in A\setminus \{0\}$. Hence $spec(G)= \{(k-1)a : a \in A\setminus \{0\}\}$. Suppose  $spec(G)$ is a subgroup of $A$, then $0 \in spec(G)$, and so $o(a)$ divides $(k-1)$. By Proposition \ref{BP4001}, we get the required result.\\ Conversely, assume that $A$ has a subgroup isomorphic to $\mathbb{Z}_p$, where $p$ is a prime number and $p| k-1$.  Let $a,b \in spec(G)$, where $a=(k-1)a_1, b=(k-1)b_1$ and $a_1,b_1 \in A\setminus\{0\}$. In both cases $a_1 =-b_1$ and $a_1 \neq -b_1$, $a+b \in spec(G)$. By Theorem \ref{BP4017}, we obtain the result.      
\end{proof}  

The upcoming result is an immediate consequence of Theorem \ref{BP4005}.

\begin{corollary}
The spec$(K_n,A)$ is subgroup of $A$, where $n>2$ if and only if $A$ has a subgroup isomorphic to $\mathbb{Z}_{p}$, where $p$ is a prime number and $p|n-1$. 
\end{corollary}

\begin{proposition}\label{BP4018}
Let $G$ be a complete $k$-partite graph with each partite size greater than one and $A$ be an Abelian group with $|A|>2$. Then $spec(G,A)$ is subgroup of $A$.  If $G$ is complete bipartite, then $spec(G,A)= A$.  
\end{proposition}
\begin{proof}
 By Proposition \ref{BP4002}, Lemma \ref{BP2011} and from the proof of Theorem \ref{BP4005}, we get $spec(G,A) = \{(k-1)a : a\in A$\}. Clearly, $spec(G,A)$ is a subgroup of $A$.  Suppose $G$ is complete bipartite, then $spec(G,A) =\{a : a\in A\}=A$.            
\end{proof}

From Theorem \ref{BP4007} and Proposition  \ref{BP4018} we see that if $A$ is an Abelian group containing at least three elements, then $spec(G, A)=A$, where $G=C_{4n}$ or $G$ is complete bipartite graph, with each partite size greater than one. 

\begin{prob}
    Characterize graphs $G$ for which $spec(G,A) = A$, where $|A| \geq 3.$
\end{prob}

\section{Construction of a family of graph whose spectra contains zero} 
In this section, our primary focus lies in identifying graphs $G$ for which  $spec(G,A)$ constitutes a subgroup of $A$. As seen in Theorem \ref{BP4016} for the groups $\mathbb{Z}_p$ and $V_4$, $spec(G,A)$ forms a subgroup if and only if $0$ belongs to $spec(G,A)$. This observation serves as a motivation to construct family of graphs whose spectrum contains zero.

Let $G$ be a simple graph. In \cite{DFA, Vinoth}, the authors defined the following relation on $V(G)$. For any $u,v \in V(G)$, define $u\sim_G v$ if and only if $N_G(u) = N_G(v)$. Clearly, the relation $\sim_G$ is an equivalence relation on $V(G)$. Let $[u]$ be the equivalence class which contains $u$ and $S$ be the set of all equivalence classes of this relation $\sim_G$. Based on these equivalence classes, we define the reduced graph $H$ of $G$ as follows. The \emph{reduced graph} $H$ of $G$ is the graph with vertex set $V(H) = S$ and two distinct vertices $[u]$ and $[v]$ are adjacent in $H$ if and only if $u$ and $v$ are adjacent in $G$. Note that, if $V(H) = \{[u_1],[u_2],\ldots,[u_k]\}$, then $G$ is the $H$-join of $\langle[u_1]\rangle, \langle[u_2]\rangle,\ldots,\langle[u_k]\rangle$, that is, $G \cong H[\langle[u_1]\rangle, \langle[u_2]\rangle,\ldots,\langle[u_k]\rangle]$ ($\langle [u] \rangle$ denote the subgraph induced by $[u]$)
and each $[u_i]$ is an independent subset 
of $G$, we take $|[u_i]|=m_i$, where $m_i \in \mathbb{N}$ for each $i$. Clearly, $H$ is isomorphic to a subgraph of $G$ induced by $\{u_1,u_2,\ldots, u_k\}$. For example, consider the graph $G$ in Figure \ref{BP4Fig1}, its reduced graph $H$ is given in Figure \ref{BP4Fig2}. Also, in Figure \ref{BP4Fig2} we have drawn the graph $G$ of Figure \ref{BP4Fig1} along with its equivalence classes.  \\

\begin{figure}
\begin{center}
\unitlength 1mm 
\linethickness{0.4pt}
\ifx\plotpoint\undefined\newsavebox{\plotpoint}\fi 
\begin{picture}(86.75,90.475)(0,0)
\put(7.431,64.481){\circle*{2.012}}
\put(63.231,82.031){\circle*{2.012}}
\put(81.456,59.981){\circle*{2.012}}
\put(79.881,29.831){\circle*{2.012}}
\put(63.681,12.281){\circle*{2.012}}
\put(6.981,29.831){\circle*{2.012}}
\put(28.581,9.356){\circle*{2.012}}
\put(29.931,85.181){\circle*{2.012}}
\put(1.25,64.4){$v_1$}
\put(27.8,88.475){$v_2$}
\put(62.9,85.675){$v_3$}
\put(83.75,59.45){$v_4$}
\put(81.725,28.15){$v_5$}
\put(64.375,7.225){$v_6$}
\put(26.675,4.4){$v_7$}
\put(3.5,25.675){$v_8$}
\thicklines
\multiput(7.592,64.558)(.0360525452,.0336995074){609}{\line(1,0){.0360525452}}
\multiput(29.548,85.081)(.337474747,-.033737374){99}{\line(1,0){.337474747}}
\multiput(29.866,85.241)(-.03371940299,-.0824){670}{\line(0,-1){.0824}}
\multiput(29.866,85.081)(-.03330233,-1.74267442){43}{\line(0,-1){1.74267442}}
\multiput(30.026,85.558)(.03371840796,-.07297910448){1005}{\line(0,-1){.07297910448}}
\put(63.913,12.214){\line(0,1){0}}
\multiput(30.026,85.4)(.03373848238,-.03761924119){1476}{\line(0,-1){.03761924119}}
\multiput(29.866,85.081)(.06935087719,-.03370850202){741}{\line(1,0){.06935087719}}
\multiput(6.956,29.556)(.0361338983,-.0337067797){590}{\line(1,0){.0361338983}}
\multiput(28.594,9.192)(.03511839465,.03373511706){1495}{\line(1,0){.03511839465}}
\multiput(7.115,29.874)(14.5098,-.0318){5}{\line(1,0){14.5098}}
\multiput(6.637,30.033)(.08460090703,.0337324263){882}{\line(1,0){.08460090703}}
\put(40.625,.25){$G$}
\end{picture}
\end{center}
\caption{}\label{BP4Fig1}
\end{figure}

\begin{figure}
\begin{center}
\unitlength 1mm 
\linethickness{0.4pt}
\ifx\plotpoint\undefined\newsavebox{\plotpoint}\fi 
\begin{picture}(139.329,100.663)(0,0)
\put(52.329,37.212){\oval(15.5,27)[]}
\put(11.954,36.462){\oval(18.75,45)[]}
\put(89.954,36.462){\oval(18.75,45)[]}
\put(11.197,36.33){\circle*{2.236}}
\put(11.697,26.08){\circle*{2.236}}
\put(11.697,47.33){\circle*{2.236}}
\put(52.197,37.83){\circle*{2.236}}
\put(89.947,43.08){\circle*{2.236}}
\put(90.947,30.58){\circle*{2.236}}
\put(129.954,36.212){\oval(18.75,45)[]}
\put(129.947,42.83){\circle*{2.236}}
\put(130.947,30.33){\circle*{2.236}}
\put(11.697,83.83){\circle*{2.236}}
\put(52.947,83.58){\circle*{2.236}}
\put(89.197,83.58){\circle*{2.236}}
\put(129.947,83.58){\circle*{2.236}}
\thicklines
\put(11.579,83.462){\line(1,0){42}}
\put(53.079,83.212){\line(1,0){76.25}}
\multiput(11.344,47.484)(.1444326241,-.0337375887){282}{\line(1,0){.1444326241}}
\multiput(52.074,37.97)(-.83426531,-.03336735){49}{\line(-1,0){.83426531}}
\multiput(11.344,26.078)(.1184011628,.0337063953){344}{\line(1,0){.1184011628}}
\multiput(52.52,37.97)(.249733333,.0337){150}{\line(1,0){.249733333}}
\multiput(52.223,37.673)(.183009434,-.03365566){212}{\line(1,0){.183009434}}
\multiput(89.386,43.173)(2.9092857,-.0318571){14}{\line(1,0){2.9092857}}
\multiput(130.116,42.727)(-.1075856354,-.0336712707){362}{\line(-1,0){.1075856354}}
\multiput(91.17,30.538)(4.409889,-.033){9}{\line(1,0){4.409889}}
\multiput(130.859,30.241)(-.1072291667,.0336770833){384}{\line(-1,0){.1072291667}}
\qbezier(52.777,38.416)(93,85.167)(129.522,43.025)
\put(1,73.729){\framebox(137.179,26.693)[]{}}
\put(73.125,77.148){$H$}
\put(6.889,47.919){$v_1$}
\put(6.243,36.606){$v_3$}
\put(48.028,40.323){$v_2$}
\put(6.48,25.469){$v_6$}
\put(88.919,45.743){$v_4$}
\put(132.234,29.888){$v_7$}
\put(131.35,42.677){$v_5$}
\put(90.156,26.292){$v_8$}
\put(9.079,87.712){$[v_1]$}
\put(50.579,88.462){$[v_2]$}
\put(86.329,87.462){$[v_4]$}
\put(126.329,87.712){$[v_5]$}
\qbezier(52,37.625)(81.5,10.25)(131,30.375)
\qbezier(52.392,83.87)(94.199,100.663)(129.997,83.516)
\put(15.75,0){Equivalence class of $v_1$}
\put(63.75,11){$G$}
\thinlines
\put(28,4){\vector(1,-1){.07}}\multiput(18.18,14.68)(.0325,-.0358333){20}{\line(0,-1){.0358333}}
\multiput(19.48,13.246)(.0325,-.0358333){20}{\line(0,-1){.0358333}}
\multiput(20.78,11.813)(.0325,-.0358333){20}{\line(0,-1){.0358333}}
\multiput(22.08,10.38)(.0325,-.0358333){20}{\line(0,-1){.0358333}}
\multiput(23.38,8.946)(.0325,-.0358333){20}{\line(0,-1){.0358333}}
\multiput(24.68,7.513)(.0325,-.0358333){20}{\line(0,-1){.0358333}}
\multiput(25.98,6.08)(.0325,-.0358333){20}{\line(0,-1){.0358333}}
\multiput(27.28,4.646)(.0325,-.0358333){20}{\line(0,-1){.0358333}}
\end{picture}
\end{center} 
\caption{}\label{BP4Fig2}
\end{figure}

In 2023, Sabeel et al. proved the following.
\begin{corollary}[\cite{Sabeel}]
Any graph $G$ is an induced subgraph of an $A$-vertex magic graph $H$, where $A$ is a finite Abelian group.   
\end{corollary}

The  following result is a generalisation of the above result,
which has been proved only for finite Abelian groups and it uses $3|V(G)|$ vertices. In our construction we use at most $2|V(G)|$ vertices, the resulting graph becomes group vertex magic.

\begin{theorem}
The class of group vertex magic graphs is universal. In other words, every graph can be embedded as an induced subgraph of a group vertex magic graph. 
\end{theorem}
\begin{proof}
Let $H$ be any graph with $n$ vertices and $A$ be an Abelian group. Let $[u_1],[u_2],\\ {[u_3]}, \ldots,[u_t]$ where $t\leq n$ be the equivalence classes under the relation $\sim_H$ on $H$. Now add a new vertex in each of the equivalence class having an odd number of  vertices. Suppose, we add the vertex $u$ in the equivalence class $[u_i]$, then  join the new vertex $u$ to all the vertices in $N_H(u_i)$. The resulting graph $H'$ is an Eulerian graph. By Theorem \ref{BP4010} and  Lemma \ref{BP4003} the resulting graph $H'$ is  $A$ vertex magic with $0$ in $spec(H', A)$ for all non-trivial Abelian groups $A$ and $H$ is an induced subgraph of $H'$.
\end{proof}
We recall the following result, which is used to study the relationship between the spectrum of a graph and its reduced graph.
\begin{theorem}[\cite{Bala2}\label{BP2003}] 
If the reduced graph $H$ of $G$ is $A$-vertex magic, where $|A| > 2$ then $G$ is $A$-vertex magic.  
\end{theorem}

Now we have, the following observation.
\begin{observation}
Let $H$ be the reduced graph of $G$. Then spec$(H,A) \subseteq spec(G,A)$.
\end{observation}

In the following definition we introduce the reduced spectrum of a graph $G$, with respect to an Abelian group.

\begin{definition}
Let $G$ be a graph and let $H$ be its reduced graph. We say that $H$ is $A'$ magic if there exists a labeling $l:V(H)\to A$ such that $0$ can also be used to label the vertex $[u]\in V(H)$
whenever $|[u]|\geq 2$ and $\displaystyle \sum_{[u]\in N_H([v])} l([u]) = \mu$, for all $[v] \in V(H)$. The collection of all such $\mu$ is called the reduced spectrum of the graph $G$ with respect to the Abelian group $A$ and is denoted by $redspec(G,A)$. 
\end{definition}

The following Proposition shows that to compute the spectrum of a graph it is enough to compute its reduced spectrum whenever the Abelian group has at least three elemenets.

\begin{proposition}\label{BP4012}
Let $G$ be a simple graph. Then 
$spec(G,A)=redspec(G,A)$, where $A$ has at least three elements.
\end{proposition}
\begin{proof}
Let $a \in spec(G,A)$. 
Assume that $l$ is an $A$-vertex magic labeling of $G$ with the magic constant  $a$. Now, define $l' : V(H) \to A$ by $l'([u])= \displaystyle \sum_{\stackrel{v \in V(G)}{v \in [u]}} l(v)$. Then $w([u])=a$, for all $[u] \in V(H)$. Conversely, let $l'$ be a $A'$-vertex magic labeling of the reduced graph $H$ of $G$ with magic constant $a$, where $a \in redspec(G,A)$. Define $l : V(G) \to A\setminus \{0\}$ by if $|[v]| =1$, then $l(v)= l'([v])$ and if $|[v]| \geq 2$, then we label the vertices of the class using Lemma $\ref{BP2011}$ in such a way that sum of the labels of all vertices in this equivalence class is equal to $l'([v])$ in the reduced graph $H$ of G. Let $v \in V(G)$. Then
\begin{align*}
  w(v) & = \sum_{u \in N_G(v)}l(u)\\
 &= \displaystyle \sum_{[u]\in N_H([v])} l'([u])\\
 &= a.
 \end{align*}
 Since $v$ is arbitrary, we get the required result. 
\end{proof}

Based on the above theorem, we can find the spectrum of a graph with a large number of vertices by transforming the graph into its reduced graph, which has minimum vertices compared to the given graph, except for a few classes of graphs (like as Path, Generalised friendship graph). We demonstrate this in the following example. 

\begin{example}
Consider the graph $G$ in Figure \ref{BP4Fig3}. The reduced graph $H$ of $G$ is $C_4$, which is $A$-vertex magic and so by Theorem \ref{BP2003}, $G$ is also $A$-vertex magic. Under the equivalence relation $\sim_G$, each equivalence class has at least two elements and so by Proposition $\ref{BP4012}$, $spec(G,A) = redspec(G,A)$ for all Abelian groups $A$ with at least three elements. In Figure $4a$, we label the reduced graph using zero only. Now, using Proposition \ref{BP4012}, we label the vertices of $G$ so that $0 \in spec(G,A)$. By Figure $4b$, $redspec(G,A)=A\setminus \{0\}$ and using Proposition \ref{BP4012}, we have $spec(G,A) = A\setminus \{0\}$. Hence, $redspec(G,A)=spec(G,A)=A$.
\end{example}


\begin{figure}
\begin{center}
\unitlength 1mm 
\linethickness{0.4pt}
\ifx\plotpoint\undefined\newsavebox{\plotpoint}\fi 
\begin{picture}(79.3,78.504)(0,0)
\put(61.482,68.861){\circle*{1.298}}
\put(74.802,50.697){\circle*{1.298}}
\put(74.348,26.479){\circle*{1.298}}
\put(59.364,9.829){\circle*{1.298}}
\put(39.384,4.986){\circle*{1.298}}
\put(16.982,15.581){\circle*{1.298}}
\put(7.295,36.922){\circle*{1.298}}
\put(12.441,57.962){\circle*{1.298}}
\put(31.21,73.249){\circle*{1.298}}
\put(30.712,77.0){$v_1$}
\put(62.498,72.114){$v_2$}
\put(77.0,52.0){$v_3$}
\put(76.997,24.981){$v_4$}
\put(60.651,6.759){$v_5$}
\put(39.341,.25){$v_6$}
\put(39.341,-6){$G$}
\put(13.155,12.207){$v_7$}
\put(1.5,34.306){$v_8$}
\put(6.949,58.795){$v_9$}
\thicklines
\multiput(31.04,73.351)(.228440299,-.033552239){134}{\line(1,0){.228440299}}
\multiput(31.361,73.351)(.0336367347,-.278277551){245}{\line(0,-1){.278277551}}
\multiput(31.254,73.244)(-.041452954,-.033726477){457}{\line(-1,0){.041452954}}
\multiput(31.254,73.351)(.03371865204,-.03673981191){1276}{\line(0,-1){.03673981191}}
\multiput(31.147,73.029)(-.03372011252,-.05102953586){711}{\line(0,-1){.05102953586}}
\multiput(59.295,9.776)(.03334426,.97377049){61}{\line(0,1){.97377049}}
\multiput(59.402,10.204)(-.128980769,-.033615385){156}{\line(-1,0){.128980769}}
\multiput(59.51,9.883)(-.03373424069,.0345){1396}{\line(0,1){.0345}}
\multiput(59.402,9.883)(-.06504358655,.03372104608){803}{\line(-1,0){.06504358655}}
\multiput(7.172,36.961)(.0337058824,-.0740692042){289}{\line(0,-1){.0740692042}}
\multiput(16.913,15.555)(.1767253086,.033691358){324}{\line(1,0){.1767253086}}
\multiput(59.51,9.883)(.0336719101,.0375191011){445}{\line(0,1){.0375191011}}
\multiput(61.222,69.176)(.0337295285,-.045942928){403}{\line(0,-1){.045942928}}
\multiput(61.222,68.962)(-.03373833206,-.04077964805){1307}{\line(0,-1){.04077964805}}
\multiput(39.602,5.281)(.0337289272,.04367241379){1044}{\line(0,1){.04367241379}}
\multiput(17.019,15.663)(.0703427673,-.0336572327){318}{\line(1,0){.0703427673}}
\multiput(12.417,58.152)(.283456621,-.033721461){219}{\line(1,0){.283456621}}
\multiput(12.524,57.618)(.033503817,-.321091603){131}{\line(0,-1){.321091603}}
\multiput(74.815,50.767)(-.033,-1.8770769){13}{\line(0,-1){1.8770769}}
\multiput(7.066,37.068)(.1664594595,.0336584767){407}{\line(1,0){.1664594595}}
\end{picture}
\end{center}  
\caption{}\label{BP4Fig3}
\end{figure}

\begin{figure}
\begin{center}
\unitlength 1mm 
\linethickness{0.4pt}
\ifx\plotpoint\undefined\newsavebox{\plotpoint}\fi 
\begin{picture}(120.425,71.625)(0,0)
\put(111.987,32.578){\oval(16.875,40.5)[]}
\put(111.981,38.535){\circle*{2.012}}
\put(112.875,27.284){\circle*{2.012}}
\put(114.037,26.887){$-a$}
\put(113.244,38.397){$a$}
\put(41.788,31.713){\oval(16.875,40.5)[]}
\put(41.106,31.594){\circle*{2.012}}
\put(41.555,22.369){\circle*{2.012}}
\put(41.555,41.494){\circle*{2.012}}
\put(37.229,43.224){$a+b$}
\put(36.648,33.25){$-b$}
\put(36.861,24.819){$-a$}
\put(77.337,31.713){\oval(16.875,40.5)[]}
\put(77.331,37.669){\circle*{2.012}}
\put(78.231,26.419){\circle*{2.012}}
\put(72.406,39.245){$-b$}
\put(76.519,22.34){$b$}
\put(8.938,31.938){\oval(16.875,40.5)[]}
\put(8.931,37.894){\circle*{2.012}}
\put(9.831,26.644){\circle*{2.012}}
\put(5.606,40.29){$-a$}
\put(8.719,22.784){$a$}
\thicklines
\multiput(8.6,37.787)(.278289474,.033552632){114}{\line(1,0){.278289474}}
\put(40.325,41.612){\line(-2,-1){30.6}}
\put(9.725,26.313){\line(6,1){31.05}}
\multiput(40.775,31.488)(-.164690722,.033634021){194}{\line(-1,0){.164690722}}
\multiput(8.825,38.013)(.0683544304,-.0337025316){474}{\line(1,0){.0683544304}}
\multiput(41.225,22.038)(-.231382979,.033503546){141}{\line(-1,0){.231382979}}
\multiput(41.675,41.388)(.307842105,-.033552632){114}{\line(1,0){.307842105}}
\multiput(76.769,37.563)(-.20816092,-.03362069){174}{\line(-1,0){.20816092}}
\multiput(40.549,31.713)(.252517007,-.033680272){147}{\line(1,0){.252517007}}
\multiput(77.669,26.762)(-.271970149,-.03358209){134}{\line(-1,0){.271970149}}
\multiput(41.225,22.262)(.0797907489,.0337026432){454}{\line(1,0){.0797907489}}
\multiput(41.45,41.163)(.0834539171,-.0336981567){434}{\line(1,0){.0834539171}}
\multiput(77.224,37.337)(.1193707483,-.033670068){294}{\line(1,0){.1193707483}}
\multiput(112.319,27.438)(-1.2747778,-.0333333){27}{\line(-1,0){1.2747778}}
\multiput(77.9,26.538)(.0940969529,.0336565097){361}{\line(1,0){.0940969529}}
\multiput(111.869,38.688)(-1.2831481,-.0333704){27}{\line(-1,0){1.2831481}}
\qbezier(8.825,38.237)(45.163,63.662)(111.647,38.688)
\qbezier(9.5,26.988)(55.063,1)(112.769,27.663)
\qbezier(8.607,37.795)(45.358,63.41)(112.653,27.294)
\qbezier(9.563,26.34)(57.052,2.714)(111.861,38.59)
\put(34.249,.25){Figure $4(a)$: magic constant zero}
\thinlines
\put(8.75,68.375){\line(1,0){103.75}}
\put(9.006,68.256){\circle*{2.012}}
\put(42.63,68.381){\circle*{2.012}}
\put(75.88,68.381){\circle*{2.012}}
\put(112.505,68.381){\circle*{2.012}}
\put(7.75,62.75){$0$}
\qbezier(9,68.125)(48.499,53.5)(112.499,68.375)
\put(41.549,62.95){$0$}
\put(74.749,63.15){$0$}
\put(111.499,63.25){$0$}
\put(6.25,56){\framebox(111.375,15.625)[]{}}
\put(56.125,56.875){$H$}
\end{picture}
\end{center}  
\end{figure}

\begin{figure}
\begin{center}
\unitlength 1mm 
\linethickness{0.4pt}
\ifx\plotpoint\undefined\newsavebox{\plotpoint}\fi 
\begin{picture}(119.925,69)(0,0)
\put(111.487,32.078){\oval(16.875,40.5)[]}
\put(111.481,38.035){\circle*{2.012}}
\put(112.376,26.784){\circle*{2.012}}
\put(113.537,26.387){$a$}
\put(112.744,37.897){$b$}
\put(41.288,31.213){\oval(16.875,40.5)[]}
\put(40.605,31.094){\circle*{2.012}}
\put(41.055,21.869){\circle*{2.012}}
\put(41.055,40.994){\circle*{2.012}}
\put(39.799,43.524){$a$}
\put(36.148,33.342){$-b$}
\put(36.361,24.319){$-a$}
\put(76.837,31.213){\oval(16.875,40.5)[]}
\put(76.831,37.169){\circle*{2.012}}
\put(77.731,25.919){\circle*{2.012}}
\put(64.404,39.565){$-(a+b)$}
\put(77.019,22.059){$a$}
\put(8.438,31.438){\oval(16.875,40.5)[]}
\put(8.431,37.394){\circle*{2.012}}
\put(9.331,26.144){\circle*{2.012}}
\put(7.006,39.89){$a$}
\put(8.019,22.284){$b$}
\thicklines
\multiput(8.1,37.287)(.278289474,.033552632){114}{\line(1,0){.278289474}}
\put(39.825,41.112){\line(-2,-1){30.6}}
\put(9.225,25.813){\line(6,1){31.05}}
\multiput(40.275,30.988)(-.164690722,.033634021){194}{\line(-1,0){.164690722}}
\multiput(8.325,37.513)(.0683544304,-.0337025316){474}{\line(1,0){.0683544304}}
\multiput(40.725,21.538)(-.231382979,.033503546){141}{\line(-1,0){.231382979}}
\multiput(41.175,40.888)(.307859649,-.033552632){114}{\line(1,0){.307859649}}
\multiput(76.271,37.063)(-.208172414,-.03362069){174}{\line(-1,0){.208172414}}
\multiput(40.049,31.213)(.25252381,-.033680272){147}{\line(1,0){.25252381}}
\multiput(77.17,26.262)(-.271977612,-.03358209){134}{\line(-1,0){.271977612}}
\multiput(40.725,21.762)(.0797907489,.0337026432){454}{\line(1,0){.0797907489}}
\multiput(40.95,40.663)(.0834562212,-.0336981567){434}{\line(1,0){.0834562212}}
\multiput(76.724,36.837)(.119377551,-.033670068){294}{\line(1,0){.119377551}}
\multiput(111.821,26.938)(-1.2749259,-.0333333){27}{\line(-1,0){1.2749259}}
\multiput(77.398,26.038)(.0941052632,.0336565097){361}{\line(1,0){.0941052632}}
\multiput(111.37,38.188)(-1.2831852,-.0333704){27}{\line(-1,0){1.2831852}}
\qbezier(8.325,37.737)(44.662,63.162)(111.147,38.188)
\qbezier(9,26.488)(54.563,.5)(112.271,27.163)
\qbezier(8.107,37.295)(44.858,62.91)(112.155,26.794)
\qbezier(9.063,25.84)(56.552,2.214)(111.361,38.09)
\put(4.5,3.5){Figure $4(b)$: magic constant $a \neq 0$ and $a,b \in A\setminus \{0\}$ $\&$ $a \neq -b$.}
\thinlines
\put(8.25,65.75){\line(1,0){103.75}}

\put(8.506,65.631){\circle*{2.012}}

\put(42.13,65.756){\circle*{2.012}}
\put(75.38,65.756){\circle*{2.012}}
\put(112.005,65.756){\circle*{2.012}}
\put(6.25,60.125){$a+b$}
\qbezier(8.5,65.5)(47.999,50.875)(111.999,65.75)
\put(38.749,60.625){$-b$}
\put(72.0,60.625){$-b$}
\put(106.999,60.625){$a+b$}
\put(5.75,53.375){\framebox(111.375,15.625)[]{}}
\put(55.625,54.25){$H$}
\end{picture}
\end{center}
\end{figure}

In the upcoming Lemma a sufficient condition for a graph $G$ to be $A$-vertex magic is given, where $A$ has at least three elements.

\begin{lemma}\label{BP4003}
Let $G$ be a graph and $\sim_G$ be  the equivalence relation on $G$. If each equivalence class has at least two elements, then $0 \in spec(G,A)$ for all $|A|>2$.
\end{lemma}
\begin{proof}
By Lemma $\ref{BP2011}$, we can label the vertices in each equivalence whose sum  equal to zero. Hence $0 \in spec(G,A)$ for all $|A|>2$.
\end{proof}

\begin{remark}
Let $G$ be a graph and $\sim_G$ be the equivalence relation on $G$ and let $A = \mathbb{Z}_2 \times \mathbb{Z}_2$ or $\mathbb{Z}_p$, where $p > 2$ is a prime number. If each equivalence class has at least two elements, then by Lemma \ref{BP4003} and Theorem \ref{BP4016}, $spec(G,A)$ is a subgroup of $A$.
\end{remark}

The following Lemma provides an infinite number of graphs whose spectra contain zero.

\begin{lemma}\label{BP4004}
Let $H$ be a simple graph on $k$ vertices and $G = H[{K_{n_1}}^c,{K_{n_2}}^c, \ldots,{K_{n_k}}^c]$. If $n_j >1$, for $1\leq j\leq k$, then $0 \in spec(G,A)$, where $|A|>2$ .
\end{lemma}
\begin{proof}
Under the equivalence relation $\sim_G$ on $G$ each equivalence class has at least two elements, now by Lemma $\ref{BP4003}$, we  get a labeling such that $w(v)=0$, for all $v \in V(G)$.
\end{proof}

\begin{theorem}\label{BP4011}
Let $P_k$ be the path on $k$ vertices and $G = P_k[{K_{n_1}}^c,{K_{n_2}}^c, \ldots,{K_{n_k}}^c]$ be a  graph, where  $k$ is even. Then $0 \in spec(G,A)$, where $|A|>2$ if and only if $n_j >1$, for each $j=1,2,\ldots,k$.
\end{theorem}
\begin{proof}
For $1 \leq j \leq k$, let $[u_j]$ denote the equivalence classes of $G$ under the equivalence relation $\sim_G$ on $G$  with $|[u_j]|=n_j$. Assume that $0 \in spec(G,A)$. Let $v \in [u_k]$. Since $w(v)=\displaystyle \sum_{u \in [u_{k-1}]}l(u) = 0$, we have $n_{k-1} > 1$. Now, $v \in [u_{k-2}]$. Since $w(v)=\displaystyle \sum_{u \in [u_{k-1}]} l(u)+ \sum_{u \in [u_{k-3}]} l(u)= 0$, we have $n_{k-3}>1$. Continuing this argument, we get $n_i > 1$, where $i$ is odd and $1 \leq i \leq k$. Let $v \in [u_1]$. Since $w(v)= \displaystyle \sum_{u \in [u_{2}]} l(u)=0$, which implies $n_2 > 1$. Continuing this argument, we get $n_i > 1$, where $i$ is even and $1 \leq i \leq k$. The  converse part follows from Lemma \ref{BP4004}.
\end{proof}

\begin{theorem}
Let $P_k$ be the path on $k$ vertices. Let $G = P_k[{K_{n_1}}^c,{K_{n_2}}^c, \ldots,{K_{n_k}}^c]$ be a  graph, where $k$ is odd. Then $0 \in spec(G,A)$, where $|A|>2$ if and only if $n_i >1$, where $i$ is even and $1< i < k$.    
\end{theorem}
\begin{proof}
Proceeding as in the proof of Theorem $\ref{BP4011}$, we get $n_i > 1$, where $i$ is even and $1< i < k$. \\ 
Conversely, assume that $n_i >1$, where $i$ is even and $1< i < k$. Using Lemma \ref{BP2011}, we define $l : V(G) \to A\setminus \{0\}$ in such a way that  \[\sum_{u\in[u_{i}]}l(u) = \begin{cases} \text{$-a$} &\quad\text{if $i \equiv 1 (\textrm{mod}\ 4)$}\\ \text{$a$} &\quad\text{if $i \equiv 3 (\textrm{mod}\ 4)$}\\ 
\text{$0$} &\quad\text{otherwise}
\end{cases} \] Thus, $w(v)= 0$, for all $v \in V(G)$.  
\end{proof}

 Let $G$ be a graph with $k$ equivalence classes $[v_1],[v_2],\ldots,[v_k]$ under the equivalence relation $\sim_G$. Now, add any number of vertices in any equivalence class. Without loss of generality add $r_i \geq 0$ vertices in $[v_i]$, which are adjacent to each vertex in $N_G(v_i)$, $i = 1,2,\ldots , k$ the resulting graph is called $G'$.
In the following theorem, we construct infinite number of  graphs, whose spectrum is same.
\begin{theorem}
Let $A$ be an Abelian group containing at least three elements. If $G$ has at least two elements in each equivalence class, then
the graphs $G$ and $G'$ have same spectrum.
\end{theorem}
\begin{proof}
 Assume that $G$ has at least two elements in each equivalence class. Let $a \in spec(G',A)$. Using Lemma \ref{BP2011}, we construct a $A$-vertex magic labeling on $G$ with sum of all label in each equivalence class is equal to sum of all label in corresponding equivalence class in $G'$. Therefore $a \in spec(G,A)$, which implies $spec(G',A) \subseteq spec(G,A)$. By a similar argument $spec(G,A) \subseteq spec(G',A)$. Hence $spec(G,A)=spec(G',A)$. 
\end{proof}


\begin{remark}
The converse of above theorem is not true. Consider the complete graph $K_2$. Clearly, $K_2$ has two equivalence classes $[v_1]$ and $[v_2]$. Now, we add  $r_i= i-1$ vertices in $[v_i]$, for $i = 1,2$. Then the resulting graph $G'= P_3$ and hence $K_2$ and $P_3$ have same spectrum, for all $|A| > 2$.  
\end{remark}

\begin{theorem}\label{BP4006}
Let $G_k$ be a collection of simple graphs, where $k = 1,2,\ldots, n$. If  $0 \in spec(G_k,A)$, for some $k$, then $0 \in spec(\otimes_{k=1}^n G_k, A)$.    
\end{theorem}
\begin{proof}
Without loss of generality, we assume that $0 \in spec(G_1,A)$  corresponding to the magic labeling  $l$. Let $v_i \in V(G_i)$, where $i = 1,2,\ldots,n$. Assume that $N_{G_i}(v_i)=\{u_{i_1},u_{i_2}\ldots,u_{i_{n_i}}\}$. 
Then $N_{\otimes_{k=1}^n G_k}((v_1,v_2,\ldots,v_n))= \{(u_{1_{j_1}},u_{2_{j_2}},\ldots, u_{k_{j_k}},\ldots, u_{n_{j_n}}) : j_k= 1,2,\ldots,n_{k} \text{ and } k = 1,2,\ldots,n \}$. 
Now, define $l' : V(\otimes_{k=1}^n G_k) \to A\setminus\{0\}$ by \\ $l'((v_1,v_2,\ldots,v_n))=l(v_1)$. Then 
\begin{align*}
w((v_1,v_2,\ldots,v_n))&= \sum_{j_1=1}^{n_1}\cdots \sum_{j_n=1}^{n_n} l'((u_{1_{j_1}},u_{2_{j_2}},\ldots, u_{n_{j_n}}))\\ 
&= n_2\cdots n_{i} \cdots n_n\sum_{j_1=1}^{n_1}l(v_{j_i}) \\
&= 0.
\end{align*}
Since $(v_1,v_2,\ldots,v_n)$ is arbitrary, $w(v)= 0$, for all $v \in V(\otimes_{i=1}^n G_i)$. 
\end{proof}

\noindent{ We recall the following Theorem in \cite{Bala1}}.

\begin{theorem}\label{BP1113}\cite{Bala1}
The graph $P_n \otimes C_m$ is group vertex magic if and only if \item[$(i)$] $n \le 3$  (or) \item[$(ii)$] $n>3$ and $m \equiv 0 (mod \, 4)$.
\end{theorem}

\begin{remark}
Since $G_1 \otimes G_2 \cong G_2 \otimes G_1$, therefore by using Theorem \ref{BP4007} and  Theorem \ref{BP4006} the second  part of the converse of Theorem \ref{BP1113} is an immediate consequence.
\end{remark}

\begin{definition}\cite{Saravanan}
 Given simple graphs $H, G_1, G_2,\ldots,G_k$, where $k = |V(H)|$, the generalized corona product denoted by $H \tilde{\circ} \wedge_{i=1}^{k} G_i$, is the graph obtained by taking one copy of graphs $H, G_1, G_2,\ldots,G_k$ and joining the $i^{th}$ vertex of $H$ to every vertex of $G_i$. In particular, if  $G_i \cong G$, for $1 \leq i \leq k$, the graph $H \tilde{\circ} \wedge_{i=1}^{k} G_i$ is called simply corona  of $H$ and $G$, denoted by $H \circ G$.
\end{definition}

In our investigation, we have the following natural question: Does there exist a graph whose spectrum contains all the elements of the abelian group A except zero?. The following result answers the question affirmatively by providing infinite class of such graphs.

\begin{theorem}
Let $H$ be a graph of order $k$.
Let $G = H \tilde{\circ} \wedge_{i=1}^{k} K_{n_i}^c$ be a graph, where $n_i>1$, for all $i$. Then  $spec(G,A)= A\setminus \{0\}$, where $|A|>2$.
\end{theorem}
\begin{proof}
Let $V(G)= V(H) \cup  \big(\bigcup_{i=1}^kV(K_{n_i}^c)\big)$, where $V(H)=\{v_1,v_2,\ldots,v_k\}$ and $V(K_{n_i}^c)=\{v_1^i,v_2^i,\ldots,v_{n_i}^i\}$, for all $i$. Let $a \in A\setminus \{0\}$. Define $l : V(G) \to A\setminus \{0\}$ by $l(v_i)=a$ and using Lemma \ref{BP2011}, we label the vertices $v_j^i$ with the condition $\sum_{j=1}^{n_i}l(v_j^i)= - (deg_H(v_i)-1)a$. Thus $w(v)= a$ for all $v \in V(G)$. Since $a$ is arbitrary and from observation \ref{BP4013}, we get the required result.
\end{proof}

\section{Conclusion and scope} In this paper, we have introduced the concept of  $A$- vertex magic spectrum of a graph. Furthermore, we have constructed an infinite family of graphs with spectra containing zero. Additionally, we have established sufficient conditions for the presence of $0$ in $spec(G,A)$, specifically in certain product graphs. We propose the following problems for further investigation.\\
1) Identify  family of graphs whose spectrum forms a subgroup.\\
2) Establish necessary conditions for a graph's spectrum to be a subgroup. \\ 
3) Characterize graph $G$ for which $spec(G,A) = A$, where $|A| \geq 3.$ 
\section*{Acknowledgements} 
The work of S. Balamoorthy is supported by a Junior Research Fellowship from CSIR-UGC, India(UGC-Ref.No.:1085/(CSIR-UGC NET JUNE 2019)).


\end{document}